\documentclass[10pt,journal]{IEEEtran}
\ifCLASSOPTIONcompsoc
   \usepackage[nocompress]{cite}
\else
   \usepackage{cite}
\fi
\ifCLASSINFOpdf
\else
\fi
\usepackage[cmex10]{amsmath}
\usepackage{algorithmic}
\ifCLASSOPTIONcompsoc
\usepackage[tight,normalsize,sf,SF]{subfigure}
\else
\usepackage[tight,footnotesize]{subfigure}
\fi
\usepackage{url}

\usepackage{amsthm,amssymb,amscd,fancybox,ifthen,float,epsfig,subfigure}
\usepackage{graphicx}
\usepackage{setspace}
\usepackage[all]{xy}
\usepackage{comment}
\usepackage{url}

\usepackage[table]{xcolor}

\newtheorem{theorem}{Theorem}[section]

\theoremstyle{definition}

\theoremstyle{remark}
\newtheorem{remark}[theorem]{Remark}

\newcommand{\bb}{\boldsymbol b}
\newcommand{\bk}{\boldsymbol k}
\newcommand{\bm}{\boldsymbol m}
\newcommand{\bx}{\boldsymbol x}
\newcommand{\br}{\boldsymbol r}
\newcommand{\by}{\boldsymbol y}

\newcommand{\btheta}{\boldsymbol \theta}

\newcommand{\cN}{\mathcal N}
\newcommand{\cGP}{\mathcal{GP}}

\begin{document}
%
\title{Fast Direct Methods for Gaussian Processes}
%
%
%
%
\author{Sivaram Ambikasaran, Daniel Foreman-Mackey,
Leslie Greengard,~\IEEEmembership{Member,~IEEE,}\\
  David W. Hogg, and Michael O'Neil,~\IEEEmembership{Member,~IEEE}
\IEEEcompsocitemizethanks{ \IEEEcompsocthanksitem S. Ambikasaran is
  with the Courant Institute at New York University, New York, NY,
  e-mail: sivaram@cims.nyu.edu.  \IEEEcompsocthanksitem
  D. Foreman-Mackey is with the Department of Physics at New York
  University, New York, NY, e-mail: danfm@nyu.edu.
  \IEEEcompsocthanksitem D. W. Hogg is with the Department of Physics
  at New York University, New York, NY, e-mail: david.hogg@nyu.edu.
  \IEEEcompsocthanksitem L. Greengard is with the Courant Institute at
  New York University, New York, NY and the Simons Center for Data
  Analysis at the Simons Foundation, New York, NY, e-mail:
  greengard@cims.nyu.edu.  \IEEEcompsocthanksitem M. O'Neil is with
  the Courant Institute and the Polytechnic School of Engineering
  at New York University, New York, NY, e-mail:
  oneil@cims.nyu.edu.}  \thanks{Research by S. Ambikasaran and
  M. O'Neil is supported in part by the Air Force Office of Scientific
  Research under NSSEFF Program Award FA9550-10-1-0180 and
  AIG-NYU Award \#A15-0098-001. Research by
  L. Greengard is supported in part by the Simons Foundation and the
  Air Force Office of Scientific Research under NSSEFF Program Award
  FA9550-10-1-0180. Research by D. Foreman-Mackey and D. W. Hogg is
  supported in part by NASA grant NNX12AI50G, NSF grant IIS-1124794,
  the Gordon and Betty Moore Foundation, and the Alfred P. Sloan
  Foundation.}}

%
%

\markboth{}
{ \MakeLowercase{\textit{et al.}}: Bare Demo of IEEEtran.cls for Computer Society Journals}
%


\IEEEcompsoctitleabstractindextext{%
\begin{abstract}
A number of problems in probability and statistics can be addressed
using the multivariate normal (Gaussian) distribution.  In the
one-dimensional case, computing the probability for a given mean and
variance simply requires the evaluation of the corresponding Gaussian
density.  In the $n$-dimensional setting, however, it requires the
inversion of an $n \times n$ covariance matrix, $C$, as well as the
evaluation of its determinant, $\det(C)$.  In many cases, such as
regression using Gaussian processes, the covariance matrix is of the
form $C = \sigma^2 I + K$, where $K$ is computed using a specified
covariance kernel which depends on the data and additional parameters
(hyperparameters).  The matrix $C$ is typically dense, causing
standard direct methods for inversion and determinant evaluation to
require $\mathcal O(n^3)$ work. This cost is prohibitive for
large-scale modeling.  Here, we show that for the most commonly used
covariance functions, the matrix $C$ can be hierarchically factored
into a product of block low-rank updates of the identity matrix,
yielding an $\mathcal O (n\log^2 n) $ algorithm for inversion.
More importantly, we show
that this factorization enables the evaluation of the determinant
$\det(C)$, permitting the direct calculation of probabilities in high
dimensions under fairly broad assumptions on the kernel defining $K$.
Our fast algorithm brings many problems in marginalization and the
adaptation of hyperparameters within practical reach using a single
CPU core. The combination of nearly optimal scaling in terms of
problem size with high-performance computing resources will permit the
modeling of previously intractable problems.  We illustrate the
performance of the scheme on standard covariance kernels.
\end{abstract}

}

\maketitle

\IEEEdisplaynotcompsoctitleabstractindextext

%
\IEEEpeerreviewmaketitle

\section{Introduction}
%
%

%
%
%
%
\IEEEPARstart{A} common task in probability and statistics is the
computation of the numerical value of the posterior probability of
some parameters $\btheta$ conditional on some data $\bx, \by \in
\mathbb R^n$ using a multivariate Gaussian distribution. This requires the
evaluation of
\begin{equation}
p(\btheta | \bx, \by) \propto \frac{1}{ | \det(C(\bx;
\btheta)) |^{1/2}}
e^{-\frac{1}{2} \by^t C^{-1}(\bx;\btheta) \by} \, p(\btheta),
\end{equation}
where $C(\bx;\btheta)$ is an $n \times n$ symmetric, positive-definite {\it
  covariance} matrix. The explicit dependence of $C$ on particular
  parameters $\btheta$ is shown here, and may be dropped in the
  proceeding discussion.  In the one-dimensional case, $C$ is simply
the scalar variance. Thus, computing the probability requires only the
evaluation of the corresponding Gaussian.  In the $n$-dimensional
setting, however, $C$ is typically dense, so that its inversion
requires $O(n^3)$ work as does the evaluation of its determinant
$\det(C)$.  This cost is prohibitive for large $n$.

In many cases, the covariance matrix $C$ is assumed to be of the form
$C(\bx) = \sigma^2 I + K(\bx)$, where $K_{ij}(\bx) = k(x_i,x_j)$.
This happens when the model for the data assumes some sort of
uncorrelated additive measurement noise having variance $\sigma^2$
in addition to some structured covariance described by the kernel $k$.
The
function $k(x_i,x_j)$ is called the covariance function or covariance
kernel, which, in turn, can depend on additional parameters,
$\btheta$.  Covariance matrices of this form universally appear in
regression and classification problems when using Gaussian process
priors~\cite{rasmussen2006gaussian}.  Because many covariance kernels
are similar to those that arise in computational physics, a
substantial body of work over the past decades has produced a host of
relevant fast algorithms, first for the rapid application of matrices
such as $K$~\cite{greengard1987fast, fastgauss, fong2009black,
  gimbutas2003generalized, ying2004kernel}, and more recently on their
inversion~\cite{greengard2009fast, ambikasaran2013fastdirect,
  amirhossein2014fast, chandrasekaran2006fast, fong2009black,
  borm2003hierarchical, martinsson2005fast, ho2012fast}.  We do not
seek to further review the literature here, except to note that it is
still a very active area of research.

Using the approach outlined in~\cite{ambikasaran2013fastdirect}, we
will show that under suitable conditions, the matrix $C$ can be
hierarchically factored into a product of block low-rank updates of the
identity matrix, yielding an $\mathcal O (n\log^2 n) $ algorithm for
inversion.  More importantly (and perhaps somewhat surprising), we
show that our factorization enables the evaluation of the determinant,
$\det(C)$, in $\mathcal O (n\log n)$ operations.  Together, these
permit the efficient direct calculation of probabilities in high
dimensions.
Previously existing methods for inversion and determinant
evaluation were based on either rough approximation methods or
iterative methods
\cite{murray,smola,snelsongh,anitescu,anitescu2}. These schemes are
particularly ill-suited for computing determinants. Although bounds
exist for sufficiently random and diagonally dominant matrices, they
are often inadequate in the general case \cite{brent}.
We briefly review existing accelerated methods for Gaussian processes
in Section~\ref{sec-accel} and present a cursory heuristic comparison
with our covariance matrix factorization.

Gaussian processes are the tool of choice for many statistical
inference or decision theory problems in machine learning and the
physical sciences.  They are ideal when requirements include
flexibility for the modeling of continuous functions. However,
applications are limited by the computational cost of matrix inversion
and determinant calculation.
Furthermore, the determinant of the covariance matrix is 
required for Gaussian
process likelihood evaluations (i.e., computation of any actual value
of the probability of the data under the covariance hyperparameters,
or evidence).
Existing linear algebraic schemes for direct matrix inversion and
determinant calculation are prohibitively expensive when the
likelihood evaluation is placed {\it inside} an outer optimization or
Markov chain Monte Carlo (MCMC) sampling loop.

In this paper we will focus on describing and applying our new
methods for handling large-scale covariance matrices (dense and full-
or high-rank) to avoid the computational bottlenecks encounters in
regression, classification, and other problems when using Gaussian
process models.
We motivate the algorithms by explaining where their need arises only 
in Gaussian process regression, but similar calculations are
frequently encountered in other regimes under Gaussian process priors.
Other applications, such as marginalization and adaptation of
hyperparameters are relatively straightforward, and the computational
bottlenecks of each are highly related.

The paper is organized as follows. Section~\ref{sec-gp} reviews some
basic facts about Gaussian processes and the resulting formulas
encountered in the case of a one-dimensional regression
problem. Prediction, marginalization, adaptation of hyperparameters,
and existing approximate accelerated methods are also discussed.
Section~\ref{sec-matrices} discusses the newly developed matrix
factorization for Hierarchical Off-Diagonal Low-Rank (HODLR) matrices,
for which factorization requires only $\mathcal O (n \log^2 n)$
work. Subsequent applications of the operator and its inverse scale as
$\mathcal O (n \log n)$.  Many popular covariance functions used for
Gaussian processes yield covariance matrices satisfying the HODLR
requirements. While other hierarchical methods could be used for this
step, we focus on the HODLR decomposition because of its simplicity
and applicability to a wide range of covariance functions. We would
like to emphasize that the algorithm will work for any covariance
kernel, but the scaling of the algorithm might not be optimal; for
instance, if the covariance kernel has a singularity or is highly
oscillatory without damping. Further, in Section~\ref{sec-deter}, we
show that the determinant of an HODLR decomposition can be computed in
$\mathcal O(n \log n)$ operations. Section~\ref{sec-numerical}
contains numerical results for our method applied to some standard
covariance functions for data embedded in varying 
dimensions.
Finally, in Section~\ref{sec-conclusions}, we summarize our results
and discuss any shortcomings and other applications of the method, as
well as future avenues of research.

The conclusion contains a cursory description of the corresponding software
packages in C++ and Python which implement the numerical schemes of
this work. These open-source software packages have been made
available since the time of submission.

\section{Gaussian processes and regression}
\label{sec-gp}

In the past two decades, Gaussian processes have gained popularity in
the fields of machine learning and data analysis for their flexibility
and robustness. Often cited as a competitive alternative to neural
networks because of their rich mathematical and statistical
underpinnings, practical use in large-scale problems remains out of
reach due to computational complexity.  Existing direct computational
methods for manipulations involving large-scale covariance matrices
require $\mathcal O (n^3)$ calculations. This causes
regression/prediction, parameter marginalization, and optimization of
hyperparameters to be intractable problems.  This scaling can be
reduced in special cases via several approximation methods, discussed
in Section~\ref{sec-accel}, however for dense, highly coupled
covariance matrices no suitable direct methods have been proposed.
Here, by {\it direct method} we mean one that constructs the inverse
and determinant of a covariance matrix to within some pre-specified
numerical tolerance {\it directly} instead of iteratively. 
The matrix inverse can then be stored for use later, much as standard
$LU$ or $QR$ factorizations.
The numerical tolerance can be measured in the spectral or Frobenius
norms, and our algorithm is able to easily achieve approximations on
the order of $10^{-12}$. Often, near machine precision ($\sim
10^{-15}$) is attainable.

The following sections contain an overview of regression via Gaussian
processes and the large computational tasks that are required at each
step.  The one-dimensional regression case is discussed for
simplicity, but similar formulae for higher dimensions and
classification problems are straightforward to derive.  In higher
dimensions, the corresponding computational methods scale with the
same asymptotic complexity, albeit with larger constants.  For a
thorough treatment of regression using Gaussian processes, see
\cite{rasmussen2006gaussian, mackay1998introduction}.

The canonical linear regression problem we will analyze assumes a model
of the form
\begin{equation}
y = f(x) + \epsilon,
\end{equation}
where $\epsilon \sim \cN(0,\sigma^2_\epsilon)$ is some form of
uncorrelated measurement noise.
Given a dataset $\{x_i, y_i \}$, the goal is to infer $f$, or
equivalently some set of parameters that $f$ depends on.
We will enforce the prior distribution of the unknown function $f$ to
be a Gaussian process,
\begin{equation}
f \sim \cGP(m, k)
\end{equation}
where $k=k(x,x')$ is some admissable covariance function ($k$
corresponds to positive definite covariance matrices), possibly
depending on some unknown hyperparameters,
and $m=m(x)$ is the expected mean of $f$.
The task of fitting hyperparameters is discussed in
Sections~\ref{sec-hyp} and~\ref{sec-adap} of the paper.
Table~\ref{table_kernel_functions} lists some of the frequently used
covariance functions for Gaussian processes.

\begin{table}[!t]
  \caption{\footnotesize Common covariance kernels used in 
    Gaussian processes}
  \begin{center}
    \rowcolors{1}{gray!30}{white}
    \resizebox{8cm}{!}{
      \begin{tabular}{|c|c|} \hline
        \textbf{Name} & \textbf{Covariance function}\\ \hline
        \hline
        \textbf{Ornstein-Uhlenbeck} & $\exp\left(-\vert x-y \vert\right)$\\
        \hline
        \textbf{Gaussian} & $\exp\left(-\vert x-y \vert^2 \right)$ \\
        \hline
        \textbf{Mat\'{e}rn family} &
        $\sigma^2
        \dfrac{ \left( \sqrt{2\nu} | x-y | \right)^{\nu} }
           {\Gamma(\nu) 2^{\nu - 1}}
            K_{\nu}\left(\sqrt{2\nu} \vert
            x-y \vert \right)$
            \\
        \hline
        \textbf{Rational Quadratic} & $\dfrac1{\left(1 + \vert x-y \vert^2\right)^{\alpha}}$\\
        \hline
      \end{tabular}
    }
  \end{center}
  \label{table_kernel_functions}
\end{table}

\subsection{Prediction}

One of the main uses for the previous model (especially in machine
learning) is to predict, with some estimated confidence, $f(\tilde
x)$ for some new input data point $\tilde x$.
This is equivalent to calculating the conditional distribution $\tilde
y | \bx, \by, \tilde x$.  
We will not assume any parametric form of $f$, and
enforce structure only through the observed data and the choice of the
mean function $m$ and covariance function $k$. 
Additionally, for the time being, assume that $k$ is fixed (i.e.,
hyperparameters are either fixed or absent).
Given the data $\bx = (x_1 \ x_2 \ \ldots
\ x_n)^t$ and $\by = (y_1 \ y_2 \ \ldots \ y_n)^t$, it is easy to show
that the conditional distribution (likelihood) of $\by$ is given by
\begin{equation}
\by | \bx \sim \cN\left( \bm(\bx), \,
\sigma^2_\epsilon I + K(\bx) \right),
\end{equation}
where the mean vector and covariance matrix are:
\begin{equation}
\begin{aligned}
\bm(\bx) &= \left( m(x_1) \ m(x_2) \ \ldots \ m(x_n) \right)^t , \\
K_{ij}(\bx)  &= k(x_i, x_j).
\end{aligned}
\end{equation}
The conditional distribution of a predicted function value, 
\mbox{$\tilde y=f(\tilde x)$}, can then be calculated as
\begin{equation}\label{eq_reg1}
\tilde y | \bx, \by, \tilde x \sim
\cN( \tilde f, \tilde \sigma^2 ),
\end{equation}
with
\begin{equation}\label{eq_reg2}
\begin{aligned}
\tilde f &= \bk(\tilde x, \bx )
\left(\sigma^2_\epsilon I + K(\bx) \right)^{-1} \by, \\
\tilde \sigma^2 &= k(\tilde x, \tilde x) -
\bk(\tilde x, \bx )
\left(\sigma^2_\epsilon I + K(\bx) \right)^{-1} \bk(\bx, \tilde x).
\end{aligned}
\end{equation}
In the previous formulas, for the sake of simplicity, we have assumed
that the mean function $\bm = \boldsymbol 0$.
The vector $\bk(\tilde x, \bx)$ is the column vector of covariances
between $\tilde x$ and all the known data points $\bx$, and
$\bk(\bx, \tilde x) = \bk(\tilde x, \bx)^t$
\cite{rasmussen2006gaussian}. 
We have therefore reduced the problem of prediction and confidence
estimation (in the expected value sense) down to matrix-vector
multiplications. For large $n$, the cost of inverting the matrix
$\sigma^2 I + K$ is expensive, with direct methods for dense
systems scaling as $\mathcal O(n^3)$. A direct algorithm for the rapid
inversion of this matrix is one of the main contributions of this
paper.

\subsection{Hyperparameters and marginalization}
\label{sec-hyp}

As mentioned earlier, often the covariance function $k$ used to model
the data $\bx$, $\by$ depends on some set of parameters, $\btheta$.
For example, in the case of a Gaussian covariance function
\begin{equation}
k(x,x';\btheta) = \beta + \frac{1}{\sqrt{2\pi} \sigma}
e^{\frac{-(x-x')^2}{2\sigma^2}},
\end{equation}
the column vector of hyperparameters is given by $\btheta = (\beta,\sigma)^t$.

Often hyperparameters correspond to some physically meaningful
quantity of the data, for example, a decay rate or some spatial scale.
In this case, these parameters are fixed once and for all according to
the specific physics or dynamics of the model.
On the other hand,
hyperparameters may be included for robustness or uncertainty
quantification and must be marginalized (integrated) away before the
final posterior distribution is calculated. In this case, for
{\it relevant} hyperparameters $\boldsymbol \theta$ and {\it nuisance}
parameters $\boldsymbol \eta$,
in order to compute the evidence
one must compute marginalization integrals of the form
\begin{equation}
p(\by,\bx|\btheta) \propto \int \frac{1}{|\det(C)|^{1/2}}
e^{-\frac{1}{2}(\bx^t C^{-1} \bx)} \, p(\boldsymbol\eta) \, d\boldsymbol\eta,
\end{equation}
where $C = C(\bx; \btheta, \boldsymbol \eta)$ is a covariance matrix corresponding to the Gaussian process
prior
and $p(\boldsymbol \eta)$ is some prior on $\boldsymbol\eta$.
If there are $r$ nuisance parameters, this is an $r$-dimensional
integral whose numerical integration requires $\mathcal O(q^r)$
quadrature nodes, where $q$ is roughly the number of quadrature nodes
needed for one-dimensional marginalization.  Unless $C^{-1}$ and
$\det(C)$ can be calculated rapidly for varying samples of the
nuisance parameters, the direct calculation of this integral is not
possible. Rapid algorithms for constructing $C^{-1}$
and $\det(C)$ would allow for the {\it direct} marginalization of
nuisance parameters, thereby {\it directly} constructing the
probability of the data, or the marginal evidence. This is in
contrast with several existing approximate Monte Carlo methods 
for computing the above integral 
  (e.g. importance sampling, MCMC, etc.), which are not
direct, and which converge with only half-order accuracy (i.e. the
numerical accuracy of the integral only decreases as $m^{-1/2}$,
where $m$ is the
number of Monte Carlo samples).
These sampling methods may decrease the number of inversions of $C$
for varying parameters, but do not completely avoid this cost.

\subsection{Adaptation of hyperparameters}
\label{sec-adap}

Alternatively, there exist situations in which the hyperparameters
$\btheta$ do not arise out of physical considerations, but rather one
would like to infer them as {\it best fit parameters}. This entails
minimizing some regression norm with respect to the parameters,
\[
\min_{\btheta} \, |\by - f(\bx)|
\]
or rather maximizing a parameter likelihood function (point estimation
using a Bayesian framework):
\[
\max_{\btheta} \, p(\theta|\bx,\by,f).
\]
In either case, some manner of non-linear optimization must be
performed because of the non-linear dependence of {\it every entry} of
the covariance matrix $C$ on the hyperparameters $\btheta$.

Regardless of the type of optimization scheme selected,
several evaluations of the evidence,
likelihood, and/or Gaussian regression must be performed -- each of
which requires evaluation of the inverse of the covariance matrix,
$C^{-1}$. In order to achieve the maximum rate of convergence of these
opimization algorithms, the full likelihood (or evidence) is required,
i.e. the numerical value of the determinant of $C$ is need.
Unless the determinant and inverse can be re-calculated and applied to the
data $\bx$ rapidly, optimizing over all possible $\btheta$'s is not a
computationally tractable problem. We skip the discussion of various
optimization procedures relevant to the adaptation of hyperparameters
in Gaussian processes~\cite{rasmussen2006gaussian} , but only point
out that virtually all of them require the re-computation of the inverse
covariance matrix $C^{-1}$.

\subsection{Accelerated methods}
\label{sec-accel}

A variety of linear-algebraic methods have been proposed to accelerate
either the inversion of $C = I + K$, the computation of its
determinant, or both.  If $K$ is of low-rank, say $p$, then it is
straightforward to compute $C^{-1}$ and $\det(C)$ using the the
Sherman-Morrison-Woodbury
formula~\cite{woodbury1950inverting,sherman1950adjustment,hager1989updating}
and the Sylvester determinant theorem~\cite{akritas1996various} in
$O(p^2 n)$ operations. However, unless an analytical form of the
low-rank property of the covariance kernel is known, some type of
dense numerical linear algebra must be performed.  Usually,
constructing general low-rank approximations to $n \times n$ matrices
requires at least $\mathcal O(p n^2)$ operations, where $p$ is the
{\it numerical rank} of the matrix. By numerical rank we loosely mean
that there are $p$ (normalized) singular values larger than some
specified precision, $\epsilon$. This definition of numerical rank is
consistent with spectral norm, and closely related to the Frobenius
norm.

Almost all of the dense matrix low-rank approximations construct
some suitable factorization (approximation) of $K$:
\begin{equation}\label{eq-lr}
K \approx Q_{n\times p} \, K^s_{p\times p} \, Q^T_{p \times n} \, ,
\end{equation}
where one can think of $Q^T$ as {\it compressing} the action of
$K$ onto a subset of points $\{x_{i1},\dots,x_{ip} \}$, $K^s$ as the
covariance kernel acting on that subset, and $Q$ as {\it
  interpolating} the result to the full set of $n$ points $\{ x_1,
x_2, \dots, x_n \}$ \cite{snelsongh,snelson_thesis,smola,quinonero}.

Iterative methods can also be applied. These are particularly
effective when there is a fast method to compute the necessary
matrix-vector products. For Gaussian covariance matrices, this can be
accomplished using the fast Gauss transform \cite{fastgauss} and its
higher-dimensional variants using $kd$-trees (see, for example,
\cite{yang2005,shen}).  Alternatively, when $k(x_i,x_j)$ is a
convolution kernel and the data are equispaced,
the Fast Fourier Transform (FFT) can be used to
accelerate the matrix vector product \cite{dietrich1997fast}. For
non-equispaced data, non-uniform Fast Fourier Transforms (NUFFTs) are
applicable \cite{dutt-rokhlin, dutt-rokhlin2, greengard-lee}.
As mentioned in the introduction, analysis-based fast algorithms can also
be used for specific kernels \cite{anitescu2} or treated using the
more general ``black-box" or ``kernel-independent" fast multipole
methods \cite{fong2009black, gimbutas2003generalized, ying2004kernel}.

In some instances, the previously described linear algebraic or
iterative methods
can be avoided all-together if an analytical decomposition of the
kernel is known. For example, much of the mathematical machinery
needed to develop
the fast Gauss transform \cite{fastgauss} relies on careful analysis 
of generating functions (or related expansions) of the Gaussian
kernel. For example,
\begin{equation}
e^{-(t-s)^2/\delta} = \sum_{j = 0}^\infty \frac{1}{j!} 
\left( \frac{s - s_0}{\sqrt{\delta}} \right) h_j 
\left( \frac{t - s_0}{\sqrt{\delta}} \right),
\end{equation}
expresses the Gaussian as a sum of separated functions in $s$, $t$,
centered about $s_0$, scaled by $\delta$, and where $h_j$ is the
$j^\text{th}$ degree Hermite function. Similar formulas, often
referred to as addition formulas or multipole expansions in the physics
literature, can be derived for other covariance kernels. As another
example, one could build a low-rank representation of covariance
matrices generated by the Mat\'{e}rn kernel using formulas of the
form:
\begin{equation}
K_\nu(w - z) = \sum_{j = -\infty}^\infty (-1)^j K_{\nu+j}(w) I_j(z),
\end{equation}
where $w$,$z$ are complex variables for which $|w|<|z|$ and $I_j$ is
the modified Bessel function of the second kind or order $j$. See
\cite{nist, watson} for a full treatment of formulas of this 
type.
Relationships such as the previous ones lead (almost) directly to fast
algorithms for the forward application of the associated covariance
matrices.
Directly building the inverse matrix (and evaluating the determinant)
is more complicated.

Before we move on, there are several other accelerated methods which
are popular in the Gaussian process community (namely greedy
approximations, sub-sampling, and the Nystr\"om method) 
\cite{rasmussen2006gaussian}. We would like
to briefly describe the accelerations that can be obtained by
interpreting Gaussian processes via a state-space model \cite{sarkka1,
sarkka2, sarkka3}.

Stochastic linear differential equations (causal, and driven by 
Gaussian noise) and
state-space models are intimately connected with Gaussian processes
and (stationary) covariance functions via the Wiener-Khintchine theorem
\cite{chatfield}. In particular, this observation allows one to
construct spectral density approximations to stationary covariance
kernels which in turn give rise to a corresponding state-space
process. This process can then be analyzed using Kalman filters and
other smoothers, which often have linear computational complexity time
for single point inference \cite{sarkka2}, \cite{ambikasaran-kalman}.
The accuracy of this inference lies in the quality of the spectral
density approximation, which is usually expressed as a rational function.
This finite-rank spectral density approximation via rational functions
can be interpreted much in the same way as approximating the Gaussian
process covariance matrix as in equation~\eqref{eq-lr} -- once this
finite-rank approximation is constructed, the resulting matrix inversion
scales as $\mathcal O(p^2n)$ by the Sherman-Morrison-Woodbury formula
(see Section~\ref{sec-matrices}). Applying state-space models to {\it
parameter} inference problems, instead of smoothing or functional
inference problems, is more subtle but several methods from signal
processing are useful. 
For a clear exposition on this topic, see \cite{sarkka1,cappe}.
Often the asymptotic computational cost of the state-space model
analysis will be similar to the algorithm of this paper because both
methods are using rank considerations to approximate the covariance
structure -- the algorithm of this paper uses a spatial-hierarchical
method, whereas fast state-space methods use a spectral approximation
of temporal data.

One last theme for increasing the scalability of Gaussian processes
to big data sets is to introduce some notion of sparsity
\cite{lawrence, quinonero}. Many of the previous
accelerated methods can be interpreted as introducing sparsity at the
covariance kernel level -- i.e. by approximating the matrix $K$ as a
finite rank operator. The resulting approximation is dense, but {\it
data-sparse}. Alternatively, one may introduce sparsity at the level
of the actual matrix $K$ by thresholding small elements away from the
diagonal. The resulting $K$ may retain high (or full) numerical rank,
but the actual matrix is sparse, thereby enabling sparse matrix
algebra to be performed which has reached a high level of
acceleration in modern computing environments. Sparsity may also be
introduced by the inclusion of data-generating latent variables
(related to the state-space interpretation of Gaussian processes),
similar to hidden Markov models
\cite{lawrence, cappe, sarkka1}.

It should be noted that all of the previous methods for accelerating
Gaussian process calculations involve some sort of approximation. 
Depending on the method, either the resulting covariance matrix is
approximated (using a low-rank factorization) or the actual covariance
kernel is approximated (using a low-rank representation, or by
approximating the actual Gaussian process by a finite-rank chain, as
in the case of the state-space models).
In each case, the analysis of the approximating Gaussian process is
different because the approximation take place at different levels in
the mathematics. Our accelerated direct method, which is described in
the next section, makes an approximation at the level of the
covariance matrix. This is akin to viewing the covariance matrix as a
continuous linear operator, and not an arbitrary data matrix.
Often this approximation is negligible as it
is near to machine precision in finite digit arithmetic.

Lastly, the evaluation of determinants is a somewhat different matter.
Most of the previously described accelerated approximations in this
section are unable to evaluate the determinant in less than
$\mathcal O(n^3)$ time since this is equivalent to constructing some
matrix factorization or all of the eigenvalues.
Taylor series approximations \cite{pacelesage} and Monte Carlo methods
have been suggested \cite{barrypace}, as well as conjugate
gradient-type methods combined with trace estimators \cite{anitescu2}.
For additional approximation methods, see the text
\cite{rasmussen2006gaussian}.
In general, however, it is difficult to
obtain accurate values for the determinant in a robust and reliable
manner.  Thus, the development of a fast, accurate, and direct method
is critical in making large-scale Gaussian process modeling useful for
for {\it exact inference} problems.

\section{Hierarchical matrices}
\label{sec-matrices}

A large class of dense matrices, for example, matrices arising out of
boundary integral equations~\cite{ying2009fast}, radial basis function
interpolation~\cite{ambikasaran2013fastdirect}, kernel density
estimation in machine learning, and covariance matrices in statistics
and Bayesian inversion~\cite{ambikasaran2013fast,
  ambikasaran2013large}, can be efficiently represented as data-sparse
hierarchical matrices. After a suitable ordering of columns and rows,
these matrices can be recursively sub-divided 
and certain sub-matrices at each level can be
well-represented by low-rank matrices. 

We refer the readers
to~\cite{hackbusch1999sparse,hackbusch2000sparse,grasedyck2003construction,hackbusch2002data,borm2003hierarchical,chandrasekaran2006fast,chandrasekaran2006fast1,ambikasaran2013thesis}
for more details on this approach.  Depending on the subdivision structure
and low-rank approximation technique, different hierarchical
decompositions exist. For instance, the fast multipole
method~\cite{greengard1987fast} accelerates the calculation of
long-range gravitational forces for $n$-body problems by
hierarchically compressing the associated matrix operator using
low-rank considerations.
The algorithm of this paper makes use of sorting data points according
to a $kd$-tree, which has the same formalism in arbitrary dimension.
The data is sorted recursively, one dimension at a time, yielding a
data structure which can be searched in at most $\mathcal O(n)$ time,
and often much faster. Once the sorting is completed, the data points
can be globally re-ordered according to, for example, a $Z$-order or
$Z$-curve. It is this ordering which generates a correspondence
between individual data points and matrix columns and rows.
Based on the particular covariance kernel and the 
data structure used (an adaptive versus a uniform
sorting), the resulting algorithm will perform slightly differently,
but with the same asymptotic scaling.

\begin{figure}[!htbp]
\resizebox{\hsize}{!}{
\includegraphics{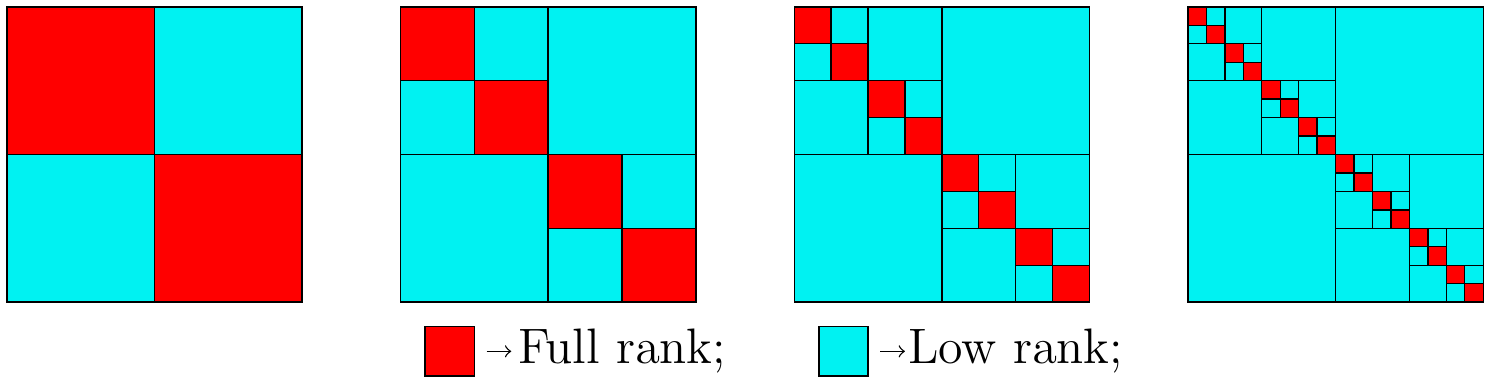}
}
\caption{ \footnotesize The same HODLR matrix at different levels.}
\label{figure_HODLR_matrices}
\end{figure}

In this article, we will be working with the class of hierarchical
matrices known as Hierarchical Off-Diagonal Low-Rank (HODLR)
matrices~\cite{ambikasaran2013fastdirect}, though the ideas extend for
other classes of hierarchical matrices as well. As the name suggests,
this class of matrices has off-diagonal blocks that are efficiently
represented in a recursive fashion. A
graphical representation of this class of matrices is shown in
Figure~\ref{figure_HODLR_matrices}. Each block represents the {\it
same} matrix, but viewed on different hierarchical scales to show the
particular rank structure.

We first give an example of a simple two-level decomposition for
real symmetric matrices, and then
describe the arbitrary-level case in more detail.
In a slight abuse of notation, in order to be consistent with previous
sources describing HODLR matrices, we will refer to the decomposition
of a matrix $K$, which is {\it not necessarily} the same $K$ as
previously mentioned in the covariance matrix case, namely in $C = I +
K$.

\begin{figure}[!htbp]
\resizebox{\hsize}{!}{
\includegraphics{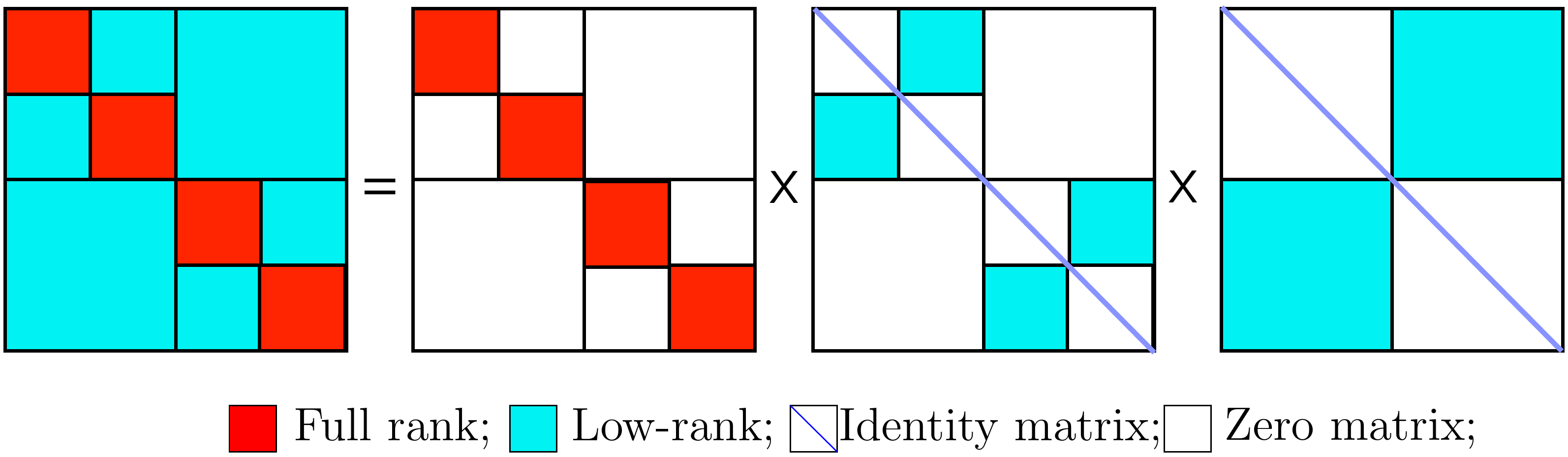}
}
\caption{ \footnotesize A two-level factorization of an HODLR matrix.}
\label{fig-two-level}
\end{figure}

Algebraically, a real symmetric matrix $K \in \mathbb{R}^{n \times n}$
is termed a
two-level HODLR matrix, if it can be written as:
\begin{equation}
K = \begin{bmatrix}
  K_1^{(1)} & U_1^{(1)} V_1^{(1)^T} \\
  V_1^{(1)} U_1^{(1)^T} & K_2^{(1)}
\end{bmatrix},
\label{equation_first_level}
\end{equation}
with the diagonal blocks given as
\begin{equation}
\begin{aligned}
K_1^{(1)} &=
\begin{bmatrix}
K_1^{(2)} & U_{1}^{(2)} 
V_1^{(2)^T}\\ V_{1}^{(2)} U_1^{(2)^T} &
K_2^{(2)}
\end{bmatrix}, \\
K_2^{(1)} &=
\begin{bmatrix}
K_3^{(2)} & U_{2}^{(2)}  V_2^{(2)^T} \\
V_{2}^{(2)} U_2^{(2)^T} & K_4^{(2)}
\end{bmatrix},
\end{aligned}
\label{equation_second_level}
\end{equation}
where the $U_i^{(j)}$, $V_i^{(j)}$ matrices are $n/2^j \times r$
matrices 
and \mbox{$r \ll n$}.  In practice, the rank of the $U$, $V$ matrices
will fluctuate slightly based on the desired accuracy of the
approximation.
In general, all off diagonal blocks of all
factors on all levels can be well-represented by a low-rank matrix,
i.e., on each level, $U_i^{(j)}, V_i^{(j)}$ are tall and thin
matrices. It is easy to show that the matrix structure given in
equations~\eqref{equation_first_level}
and~\eqref{equation_second_level} can be manipulated to provide a
factorization of the original matrix as a product of matrices, one of
which is block-diagonally dense, and the rest of which are block-diagonal 
low-rank updates to the identity matrix. This is shown in
Figure~\ref{fig-two-level}.

The above is merely a description of the structure of matrices which
meet the HODLR requirements, but not a description of how to actually
construct the factorization. There are two aspects which need to be
discussed: (i) constructing the low-rank approximations of all the
off-diagonal blocks, and (ii) using these low-rank approximations to
recursively build a factorization of the form shown in
Figure~\ref{fig-two-level}. We now describe several methods for
constructing the low-rank approximations in the next section.

\subsection{Fast low-rank approximation of off-diagonal blocks}

The first key step is to have a computationally efficient way of
obtaining the low-rank factorization of the off-diagonal blocks. Given
any matrix $A \in \mathbb{R}^{m \times n}$, the {\it optimal} low-rank
approximation (in the least-squares sense) is obtained using the
singular value decomposition (SVD)~\cite{golub1996matrix}. The
downside of using the SVD is that the computational cost of direct
factorizations scales as $\mathcal{O}(mnr)$, where $r$ is the
numerical rank of the matrix. In practice, $r$ is obtained on-the-fly
such that the factorization is accurate to some specified precision
$\epsilon$.  For our algorithm to be computationally tractable, we
need a fast low-rank factorization. More precisely, we need algorithms
that scales at most as $\mathcal{O}(r^2 n)$ to obtain a rank $r$
factorization of a $n \times n$ matrix. Thankfully, there has recently
been tremendous progress in obtaining fast low-rank factorizations of
matrices. These techniques can be broadly classified as either {\it
  analytic} or {\it linear-algebraic} techniques.

If the matrix entries are obtained as evaluations from a smooth
function, as is
the case for most of the covariance matrices in Gaussian processes, we can
rely on approximation theory based analytic techniques like 
interpolation, multipole expansion, eigenfunction expansion, Taylor
series expansions, etc. to obtain a low-rank decomposition. 
In particular, if the matrix elements are given in terms of a smooth
function $f$, as in the Gaussian process case,
\begin{equation}
A_{ij} = f(x_i,x_j),
\end{equation}
then polynomial interpolation methods can be used to efficiently
approximate the matrix $A$ with near spectral accuracy. Barycentric
interpolation formulae such as those recently discussed by Townsend
and Trefethen
and others \cite{townsend, trefethen} serve to effectively
factorize $A$ into
\begin{equation}
A \approx E \tilde A P,
\end{equation}
where $\tilde A$ is a matrix obtained by sampling the {\it function}
$f$ at suitable chosen nodes, e.g. Chebyshev interpolation nodes. The
matrices $E$, $P$ are then obtained via straightforward 
interpolation formulas. The
accuracy of the approximation can be estimated from spectral analysis
of
the interpolating Chebyshev polynomial, and the approximation can be
computed in $\mathcal O(r\max{(m,n)})$ time.

On the
other hand, if there is no {\it a-priori} information of the matrix,
then linear-algebraic methods provide an attractive way of computing fast
low-rank decompositions. These include techniques like pseudo-skeletal
approximations~\cite{goreinov1997theory}, interpolatory
decomposition~\cite{fong2009black}, randomized
algorithms~\cite{frieze2004fast,
  liberty2007randomized,woolfe2008fast}, rank-revealing
$LU$~\cite{miranian2003strong,pan2000existence}, adaptive cross
approximation~\cite{rjasanow2002adaptive,zhao-kezhong} (which is a minor
variant of partial-pivoted $LU$), and rank-revealing
$QR$~\cite{gu1996efficient}.
Though purely analytic techniques can be faster
since many operations can be pre-computed, algebraic techniques are
attractive for constructing black-box low-rank factorizations. 
The algorithm of this paper relies on an implementation of
approximate partial-pivoted $LU$, which we will now discuss.

Briefly, we construct factorizations of off-diagonal blocks via a
partial-pivoted $LU$ decomposition which executes in $\mathcal O(r
n)$ time.
Heuristically, this factorization constructs a series of rank-one 
matrices whose sum approximates the original matrix, i.e. we wish to
write
\begin{equation}
A \approx \sum_{k = 1}^r \alpha_k \boldsymbol u_k \boldsymbol v^T_k.
\end{equation}
The vectors $\boldsymbol u_k$, $\boldsymbol v_k$ are computed from the
columns and rows of $A$.

The linear complexity is achieved by checking the resulting
approximation against only a sub-sampling of the original matrix. If
the underlying matrix (covariance kernel) is sufficiently smooth, then
this sub-sampling error estimation will result in an approximation
which is accurate to near machine precision.
For other matrices or covariance kernels which
are highly oscillatory or contain small-scale structure, this method
will not scale and will likely yield a less-accurate approximation.
In this case, analytic methods are preferable as they will be more
efficient and provide suitable high-accuracy approximations.
We omit a pseudo-code description of this algorithm, as it is a
well-know linear algebra procedure, and instead refer to Section~2.2,
Algorithm~6 of \cite{rjasanow2002adaptive}.

The next section presents the fast matrix factorization of the entire
covariance matrix once the low-rank decomposition of the off-diagonal
blocks has been obtained using one of the above mentioned techniques.
We offer a concise, but complete description of the factorization in
order to make the exposition self-contained.  For a longer and more
detailed discussion of the material,
see~\cite{ambikasaran2013fastdirect}.

\subsection{HODLR matrix factorization}
\label{subsection_HODLR_matrix_factorization}

\begin{figure}[!htbp]
\resizebox{\hsize}{!}{
\includegraphics{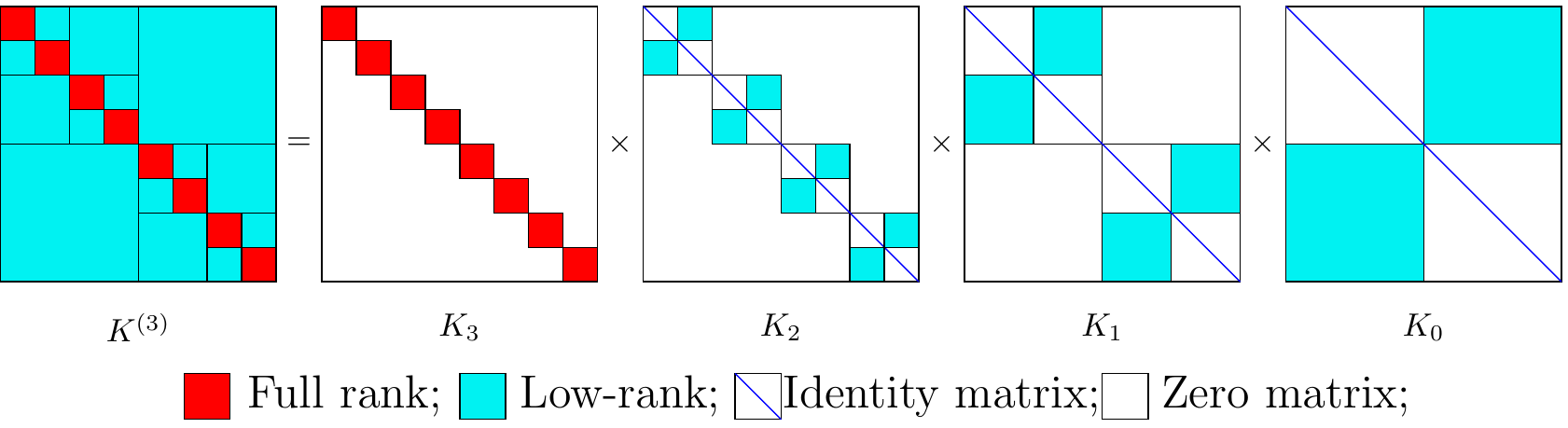}
}
\caption{ \footnotesize Factorization of a three level HODLR matrix.}
\label{figure_HODLR_level3}
\end{figure}

The overall idea behind the $\mathcal{O}(n \log^2 n)$ factorization of
an $n \times n$, $\kappa$-level (where $\kappa \sim \log n$) HODLR matrix
as described in~\cite{ambikasaran2013fastdirect} is to
factor it as a product of $\kappa+1$ block diagonal matrices,
\begin{align}
K = K_{\kappa} \, K_{\kappa-1} \, K_{\kappa-2} \cdots K_1 \, K_0,
\label{equation_HODLR_factorization}
\end{align}
where, except for $K_\kappa$, 
$K_k \in \mathbb{R}^{n \times n}$ is a block diagonal matrix
with $2^k$ diagonal blocks, each of size $n/2^k \times n/2^{k}$.
More importantly, each of these diagonal blocks is a low-rank update
to the identity matrix.
The first factor $K_\kappa$ is formed from dense block diagonal sub-matrices
of the original matrix, $K$.
Aside from straightforward block-matrix algebra,
the main tool used in constructing this factorization is the
Sherman-Morrison-Woodbury
formula~\cite{woodbury1950inverting,sherman1950adjustment,hager1989updating}.
To simplify the notation assume for a moment that $K$ is an
$n \times n$ matrix, where $n = 2^m$ for some integer $m$.
For example, a two-level HODLR matrix described in
equations~\eqref{equation_first_level}
and~\eqref{equation_second_level} can be factorized as:
\begin{align}
\label{equation_HODLR_level2_factorization}
\resizebox{\hsize}{!}{
$
\begin{bmatrix}
K_1^{(2)} & 0 & 0 & 0\\
0 & K_2^{(2)} & 0 & 0\\
0 & 0 & K_3^{(2)} & 0\\
0 & 0 & 0 & K_4^{(2)}\\
\end{bmatrix}
\begin{bmatrix}
I_{n/4} & \tilde K_{12}^{(2)} & 0 & 0 \\
\tilde K_{21}^{(2)} & I_{n/4} & 0 & 0  \\
0 & 0 & I_{n/4} & \tilde K_{34}^{(2)} \\
0 & 0 & \tilde K_{43}^{(2)} & I_{n/4}
\end{bmatrix}
\begin{bmatrix}
I_{n/2} & \tilde K_{12}^{(1)} \\
\tilde K_{21}^{(1)} & I_{n/2}
\end{bmatrix},
$
}
\end{align}
where $I_{m}$ is the $m \times m$ identity matrix, and the matrices
$\tilde{K}_{ij}^{(k)}$ are low-rank.  Similarly,
Figure~\ref{figure_HODLR_level3} graphically depicts the factorization
of a level $3$ HODLR matrix.

In the case of a one-level factorization, we can easily write down
the computation. Let the matrix $K$ be:
\begin{equation}
K = \begin{bmatrix}
A_{11} & UV^T \\
VU^T & A_{22}
\end{bmatrix},
\end{equation}
where we assume that $U$, $V$ have been computed using one of the
algorithms of the previous section. Then the {\it only} step in the
decomposition is to factor out the terms $A_{11}$, $A_{22}$, giving:
\begin{equation}\label{eq-two-level-1}
K = \begin{bmatrix}
A_{11} & 0 \\
0 & A_{22}
\end{bmatrix}
\begin{bmatrix}
I_{n/2} & A_{11}^{-1} U V^T \\
A_{22}^{-1} V U^T & I_{n/2}.
\end{bmatrix}.
\end{equation}
We see that the computation involved was to merely apply the inverse
of the dense block diagonal factor to the corresponding rows in the
remaining factor. Furthermore, since the matrix $UV^T$ was low-rank,
so is $A_{11}^{-1} UV^T$. Unfortunately, a one-level factorization
such as this is still quite expensive: it required the direct
inversion of $A_{11}$, $A_{22}$, each of which are $n/2 \times n/2$
matrices. The procedure must be done recursively across $\log n$ 
levels in order to achieve a nearly optimal algorithm.

\begin{figure*}[!t]
  \begin{equation}
\label{eq-two-level}
\resizebox{.9\hsize}{!}{$
K = \begin{bmatrix}
A_{11} &0 &0 &0 \\
0& A_{22} &0 &0 \\
0&0 & A_{33} &0 \\
0 &0 &0 & A_{44}
\end{bmatrix}
\begin{bmatrix}
I_{n/4} & A_{11}^{-1} U_1^{(2)}V_1^{{(2)}^T} &0 &0 \\
A_{22}^{-1} V_1^{(2)}U_1^{{(2)}^T} & I_{n/4}  &0 &0 \\
0 &0 & I_{n/4} & A_{33}^{-1} U_2^{(2)}V_2^{{(2)}^T} \\
0 &0 & A_{44}^{-1} V_2^{(2)}U_2^{{(2)}^T}  & I_{n/4}  
\end{bmatrix}
\begin{bmatrix}
I_{n/2} & A_1^{-1} U_1^{(1)}V_1^{{(1)}^T} \\
A_2^{-1} V_1^{(1)}U_1^{{(1)}^T} & I_{n/2}
\end{bmatrix}
$}
\end{equation}

\hrulefill
\end{figure*}

Before describing the general scheme, we give the full two-level
factorization using the notation of equations~\eqref{equation_first_level}
and~\eqref{equation_second_level}. The full factorization in this
two-level scheme is given in equation~\eqref{eq-two-level} (spanning
two columns on the proceeding page).
The matrices $A_1$ and $A_2$ appearing in the off-diagonal 
expressions are given by:
\begin{equation}\label{eq-two-level-terms}
\begin{aligned}
A_1 &= \begin{bmatrix}
A_{11} & U_1^{(2)}V_1^{{(2)}^T} \\
V_1^{(2)}U_1^{{(2)}^T} & A_{22} 
\end{bmatrix} \\
A_2 &= \begin{bmatrix}
A_{33} & U_2^{(2)}V_2^{{(2)}^T} \\
V_2^{(2)}U_2^{{(2)}^T} & A_{44} 
\end{bmatrix}
\end{aligned}
\end{equation}
This factorization is an indication of how to construct the ultimate
$\kappa$-level factorization as it only required the direct
construction of the inverse of dense matrices of size $n/4 \times
n/4$. If this procedure is repeated recursively, the only dense
inversions required are of $n/2^\kappa \times n/2^\kappa$ matrices.

At first glance, it may look as though the computation of $A_1^{-1}$,
$A_2^{-1}$ is expensive, and will scale as $\mathcal O(n^3/8)$.
However, these matrices are of the form:
\begin{equation}
\label{equation_Low_Rank_Update}
\begin{bmatrix}
A & U V^T\\
V U^T & B
\end{bmatrix}
=
\begin{bmatrix}
A & 0 \\
0 & B 
\end{bmatrix}
+ 
\begin{bmatrix}U & 0\\ 0 & V\end{bmatrix}
\begin{bmatrix}0 & V^T\\ U^T & 0\end{bmatrix}.
\end{equation}
If the inverses of $A$, $B$ are known (and they are in this case, they
were computed on a finer level), and $U$, $V$ are low-rank
matrices, then the inverse of the full matrix can be computed rapidly
using the Sherman-Morrison-Woodbury formula:
\begin{equation*}
\left( A + LSR \right)^{-1} = A^{-1} - A^{-1} L \left( S^{-1} 
+ R A^{-1} L \right)^{-1} R A^{-1}.
\end{equation*}
If $S$ is of small rank, then the inner inverse can be computed very
rapidly.

To summarize, see Figure~\ref{fig-algorithm} for rough pseudo-code
describing how to construct
a general $\kappa$-level HODLR factorization.
We avoid too much index notation, please see 
\cite{ambikasaran2013fastdirect} for a full detailed algorithm.

\begin{figure*}[!t]
\begin{algorithmic}[1]

\REQUIRE Factorization precision $\epsilon > 0$
\REQUIRE Size of smallest sub-matrix on finest level, $p_{\max}$
\COMMENT{the size of the smallest diagonal block}
\REQUIRE Matrix entry evaluation routine, $f(i,j)$
\STATE $\kappa \leftarrow \left\lfloor \log_2 n/p_{max} \right\rfloor$
\FOR{$j = 1$ to $\kappa$}
\FOR{all $i$}
\STATE Compute the low-rank factorization $U_i^{(j)}V_i^{{(j)}^T}$
of all off-diagonal blocks in the $\kappa$-level\\ 
hierarchy to precision $\epsilon$
\ENDFOR
\ENDFOR

\STATE Form the block diagonal matrix $K_\kappa$ using the block
diagonals of the original matrix $K$

\FOR {$j = \kappa-1$ down to $1$} 
\FOR {$\ell = j-1$ down to $1$}
\STATE Apply the inverse of the block diagonals of $K_j$ to each
left low-rank factor, $U_i^{(\ell)}$, $V_i^{(\ell)}$,\\ in the remaining
off-diagonal blocks \COMMENT{Factor $K_j$ out of the remaining matrix}
\ENDFOR
\ENDFOR
\STATE $K_0$ is the last factor, off the form $I + \text{low-rank}$, obtained
by having applied all earlier factor inverses\\ to off-diagonal blocks
on finer levels
\end{algorithmic}
\caption{Pseudo-code for constructing an HODLR factorization.}
\label{fig-algorithm}
\end{figure*}

This pseudo-code computes a factorization of the original
matrix $K$. We have not yet computed the inverse $K^{-1}$. The inverse
can be computed by directly applying the Sherman-Morrison-Woodbury
formula to each term in the factorization
\begin{equation}\label{eq-fac}
K = K_\kappa \, K_{\kappa - 1} \, \cdots \, K_1 \, K_0.
\end{equation}
Since each term is block diagonal or a block diagonal low-rank update
to the identity matrix, the inverse factorization
can be computed in $\mathcal O(n \log n)$ time.

Before moving on
we would like to point out that in the case where the data points at
which the kernel is to be evaluated at are not approximately uniformly
distributed, the performance of the factorization may suffer, but only
slightly. A higher level of compression could be obtained in the
off-diagonal blocks if the hierarchical tree structure is constructed
based on spatial considerations instead of {\it point count}, as is
the case with some $kd$-tree implementations.

The next section gives a brief estimate of the computational
complexity of constructing a HODLR-type factorization. 

\subsection{Computational complexity}

Constructing a HODLR-type factorization can be split into two main
steps: (i) computing the low-rank factorization of all off-diagonal
blocks, and (ii) using these low-rank approximations to recursively
factor the matrix into roughly $\mathcal O(\log n)$ pieces.

For an $n\times n $ matrix which admits the HODLR structure, as shown in, 
Figure~\ref{figure_HODLR_matrices}, there are approximately 
$\kappa \approx \log_2 n/p$, where $p$ is the size of the diagonal block
on the finest level (this is a user-defined parameter). Ignoring the
diagonal blocks, this means
there are two blocks of size $n/2 \times n/2$, four blocks of size
$n/4 \times n/4$, etc. Finding the low-rank approximation of an $n/2^j
\times n/2^j$ off-diagonal block using cross approximation requires
$\mathcal O(r n/2^j)$ flops, where $r$ is the $\epsilon$-rank of the
sub-matrix. Constructing all such factorizations requires $\mathcal
O(r n \log n)$.

Once the approximations are obtained, the matrix must be pulled apart
into its HODLR factorization. Let us remember that there are $\kappa
\approx \log_2 n$ levels in the HODLR structure. In order to factor
$K_\kappa$ (the matrix of dense block diagonals), as in
equation~\eqref{equation_HODLR_factorization}, 
out of the original matrix $K$, we must apply the
inverse of the corresponding block diagonal to all the left low-rank factors,
$U^{(\kappa)}_i$, $V^{(\kappa)}_i$ as in 
equations~\eqref{eq-two-level-1}~-~\eqref{eq-two-level-terms}.
In general, a $p \times p$ inverse must computed and applied to all
left low-rank factors, of which there are $\mathcal O(\kappa)$.
The inverse calculation is $\mathcal O(p^3)$,
and the subsequent application is $\mathcal O(prn/2^j)$, $j =
1,\ldots,\kappa$, which
dominates the inverse calculation. There
are $2^\kappa$ such applications, yielding the cost {\it for only the
first factorization level} to be $\mathcal O(prn\log n)$.
Applying the same reasoning as each subsequent factor $K_j$, $j =
\kappa-1, \ldots, 0$, along with the complexity result,
\begin{equation}
\sum_{j = 1}^{\log n} j = \mathcal O(\log^2 n),
\end{equation}
yields a total complexity for the factorization stage of $\mathcal O(n
\log^2 n)$. The factors of $p$, $r$ have been dropped as it is assumed
that $p, r \ll n$.

Given a HODLR-type factorization, it is straightforward to show that
the computational complexity of determining the inverse 
scales as $\mathcal O(n \log n)$. There
are $\mathcal O(\log n)$ factors, and since each level is constructed
as low-rank updates to the identity, invoking the
Sherman-Morrison-Woodbury formula yields the inverse in $\mathcal
O(n)$ time. This gives a total runtime of $\mathcal O(n \log n)$.

\section{Determinant computation}
\label{sec-deter}

As discussed earlier, once the HODLR factorization has been obtained,
Sylvester's determinant theorem~\cite{akritas1996various} enables the
computation of the determinant at a cost of $\mathcal{O}(n \log n)$
operations. This computationally inexpensive method for direct
determinant evaluation enables the efficient {\it direct} evaluation
of probabilities. We now briefly review the algorithm used for
determinant evaluation.

\begin{theorem}[Sylvester's Determinant Theorem]
\label{lemma_Sylvester}
If $A \in \mathbb{R}^{m \times n}$ and $B \in \mathbb{R}^{n \times
  m}$, then
$$\det\left(I_m + AB\right) = \det \left(I_n + BA\right),$$ where $I_k
\in \mathbb{R}^{k \times k}$ is the identity matrix.  In particular,
for a rank $p$ update to the identity matrix,
\[
\det(I_n+U_{n \times p}V_{p \times n}) = \det(I_p +V_{p \times n}U_{n \times p}).
\]
\end{theorem}

\begin{remark}
\label{remark_det_cost}
The computational cost associated with computing the determinant of a
rank $p$ update to the identity is $\mathcal{O}(p^2n)$.  The dominant
cost is computing the matrix-matrix product $V_{p \times n}U_{n \times
  p}$.
\end{remark}

Furthermore, we recall two basic facts regarding the determinant. First,
the determinant of a block diagonal matrix is the product
of the determinants of the individual blocks of the matrix. Second,
the determinant of a square matrix is completely multiplicative over the
set of square matrices, that is to say,
\begin{equation}
\label{remark_det_prod}
\det(A_1 A_2 \cdots A_n) = \det(A_1)
\det(A_2) \cdots \det(A_n).
\end{equation}
Using the HODLR factorization in equation~\eqref{equation_HODLR_factorization}
and these two facts, we have:
\begin{multline}
\label{equation_determinant_completely_multiplicative}
\det(K) = \det(K_{\kappa}) \det(K_{\kappa-1}) \det( K_{\kappa-2})
\cdots \\ \det(K_2) \det( K_1) \det( K_0).
\end{multline}
Each of the determinants on the right hand side of
equation~\eqref{equation_determinant_completely_multiplicative} can be
computed as a product of the determinants of diagonal blocks.  Each of
these diagonal blocks is a low-rank perturbation to the identity, and
hence, the determinant can be computed using Sylvester's Determinant
Theorem~\ref{lemma_Sylvester}. It is easy to check that the
computational cost for each of the determinants $\det(K_i)$ is
$\mathcal{O}(n)$, therefore the total computational cost for obtaining
$\det(K)$ is $\mathcal{O}(\kappa n)$, and we recall that $\kappa \sim
\log n$.

\section{Numerical results}
\label{sec-numerical}

In this section, we discuss the performance of the previously described
algorithm for the inversion and application of covariance matrices
$C = I + K$, as well as the calculation of the normalization factor, i.e., the
determinant $\det(C)$.  Detailed results for $n$-dimensional
datasets $\bx$, where each $x_i$ is a point in one, two, and three
dimensions, are provided for the Gaussian and multiquadric covariance
kernels.  We also provide benchmarks in one dimension for covariance
matrices constructed from exponential, inverse multiquadric, and
biharmonic covariance functions. Unless otherwise, stated,
all the proceeding numerical
experiments have been run on a MacBook Air with a 1.3GHz Intel Core i5
processor and 4 GB 1600 MHz DDR3 RAM. In all these cases, the matrix
entry $C_{ij}$ is given as
\[
C_{ij} = \sigma^2_i \delta_{ij} + k(r_i,r_j)
\]
where $k(r_i,r_j)$ is a particular covariance function evaluated at two
points, $r_i$ and $r_j$.  
It should be noted that in certain cases, the
matrix $K$ with entries $K_{ij}(\br) = k(r_i,r_k)$ might itself be a
rank deficient matrix and that this has been exploited in the past to
construct fast schemes. However, this is not always the case, for
instance, if the covariance function is an exponential. In this
situation, the rank of $K$ is in fact full-rank. Even if the matrix
$K$ were to be formally rank deficient, in practice, the rank might be
very large whereas the ranks of the off-diagonal blocks in the
hierarchical structure are very small. Another major advantage of
this hierarchical approach is that it is applicable to a wide range of
covariance functions and can be used in a black-box, plug-and-play,
fashion.

In each of the following tables, timings are provided for the
assembly, factorization, and inversion of the covariance matrix
(denoted by the columns {\it Assembly}, {\it Factor}, and {\it
  Solve}). Assembly of the matrix refers to computing all of the
low-rank factorizations of the off-diagonal blocks (including the time
required to construct a $kd$-tree on the data), and Factor refers to
computing the $\log n$ level factorization described earlier.  The
time for computing the determinant of the $n \times n$ matrix is given
in the {\it Det.} column, and the error provided is approximately the
relative $l_2$ precision in the solution to a test problem $C\bx =
\bb$, where $\bb$ was generated {\it a priori} from the known vector
$\bx$.

\subsection{Gaussian covariance}

\begin{figure}
  \begin{center}
    \includegraphics[width=.95\linewidth]{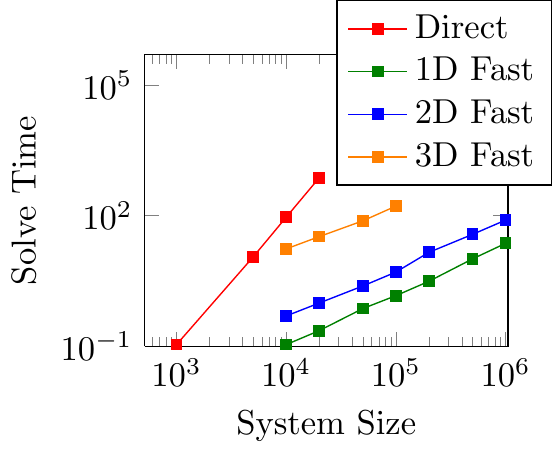}
  \end{center}
  \caption{ \footnotesize Comparison of time required to factorize the
    covariance matrix in the case of a Gaussian covariance kernel in
    one, two, and three dimensions. The conventional direct
    calculation is independent of dimension but scales as $\mathcal
    O(n^3)$.}
\label{Figure_Gaussian}
\end{figure}

Here the covariance function, and the corresponding entry of $K$,
is given as
\begin{equation}
k(r_i,r_j) = \exp \left(-{| r_i - r_j |^2}\right),
\end{equation}
where $r_i$, $r_j$ are points in one, two, or three dimensions.
The results for the Gaussian covariance kernel have been aggregated in
Table~\ref{tab_gauss_all}.
Scaling of the algorithm for data embedded in
one, two, and three dimensions is compared with the direct calculation
in Figure~\ref{Figure_Gaussian}.

\begin{table*}[!tb]
      \caption{ \footnotesize Timings for Gaussian covariance
        functions in one, two, and three dimensions. The matrix
        entries are given as $C_{ij} = 2\delta_{ij}
        + \exp \left(-|| r_i - r_j||^2\right)$, where $r_i$ are random
        uniformly distributed points in the interval $[-3, 3]^d$ ($d
        = 1, 2, 3$).}
\label{tab_gauss_all}
    \rowcolors{1}{gray!30}{white}
    \begin{center}
  \resizebox{\hsize}{!}{
    \begin{tabular}{|r|ccccc|ccccc|ccccc|}
        \hline
        & \multicolumn{5}{|c|}{One-dimensional data} &
        \multicolumn{5}{|c|}{Two-dimensional data} &
        \multicolumn{5}{|c|}{Three-dimensional data} \\
        \hline
        $n$ & Assembly & Factor & Solve & Det. & Error
        & Assembly & Factor & Solve & Det. & Error
        & Assembly & Factor & Solve & Det. & Error\\
        \hline
        $10,000$ & $0.12$ & $0.11$ &  $0.008$ & $0.01$ &$10^{-13}$ 
          & $0.56$ & $0.50$ & $0.018$ & $0.03$ & $10^{-13}$
          & $15.4$  & $17.3$ & $0.113$ & $0.91$ & $10^{-12}$ \\ 
          \hline
        $20,000$ & $0.15$ & $0.23$ &  $0.016$ & $0.03$ & $10^{-13}$ 
          & $1.16$ & $0.99$ & $0.028$ & $0.05$ & $10^{-13}$
          & $30.9$  & $33.1$ & $0.224$ & $1.06$ & $10^{-12}$ \\ 
          \hline
        $50,000$ & $0.47$ & $0.71$ &  $0.036$ & $0.12$ & $10^{-12}$
          & $2.74$ & $2.44$ & $0.067$ & $0.12$ & $10^{-13}$
          & $75.5$  & $76.3$ & $0.434$ & $1.68$ & $10^{-11}$ \\
          \hline
        $100,000$ & $1.24$ & $1.46$ &  $0.052$ & $0.24$ & $10^{-12}$
          & $5.43$ & $5.08$ & $0.165$ & $0.23$ & $10^{-12}$
          & $149$   & $166$  & $0.923$ & $3.11$ & $10^{-11}$ \\ 
          \hline
        $200,000$ & $2.14$ & $3.12$ &  $0.121$ & $0.39$ & $10^{-13}$
          & $12.4$ & $14.4$ & $0.485$ & $0.44$ & $10^{-12}$&&&&& \\
          \hline
        $500,000$ & $6.13$ & $10.2$ &  $0.388$ & $0.56$ & $10^{-12}$
          & $31.7$ & $37.3$ & $1.33$  & $1.17$ & $10^{-12}$&&&&& \\ 
          \hline
        $1,000,000$ & $14.1$ & $23.2$ &  $0.834$ & $1.52$ & $10^{-12}$
          & $70.8$ & $79.2$ & $3.15$  & $2.24$ & $10^{-12}$&&&&& \\
          \hline
      \end{tabular}
  }
    \end{center}
\end{table*}

\subsection{Multiquadric covariance matrices}

Covariance functions of the form
\begin{equation}
k(r_i,r_j) = \sqrt{\left(1+{|
    r_i - r_j |^2}\right)}
\end{equation}
are known as multiquadric covariance functions, one class of
frequently used radial basis functions. Analogous numerical results
are presented below in one, two, and three dimensions as were in the
previous section for the Gaussian covariance
function. Table~\ref{tab_multi_all} contains the results.
Scaling is virtually identical to
the Gaussian case, and we omit the corresponding plot.

\begin{table*}[!tb]
      \caption{ \footnotesize Timings for multiquadric covariance
        functions in one, two, and three dimensions. The matrix
        entries are given as $C_{ij} = \delta_{ij}
        + \sqrt{ 1 + || r_i - r_j||}$, where $r_i$ are random
        uniformly distributed points in the interval $[-3, 3]^d$ ($d
        = 1, 2, 3$).}
\label{tab_multi_all}
    \rowcolors{1}{gray!30}{white}
    \begin{center}
  \resizebox{\hsize}{!}{
    \begin{tabular}{|r|ccccc|ccccc|ccccc|}
        \hline
        & \multicolumn{5}{|c|}{One-dimensional data} &
        \multicolumn{5}{|c|}{Two-dimensional data} &
        \multicolumn{5}{|c|}{Three-dimensional data} \\
        \hline
        $n$ & Assembly & Factor & Solve & Det. & Error
        & Assembly & Factor & Solve & Det. & Error
        & Assembly & Factor & Solve & Det. & Error\\ \hline
      $10,000$    & 0.08 & 0.13 & 0.006 & 0.02 & $10^{-13}$ 
        & 0.77 & 0.86 & 0.022 & 0.04 & $10^{-13}$ 
        & 19.0 & 23.2 & 0.135 & 1.32 & $10^{-12}$ \\ \hline
      $20,000$    & 0.11 & 0.22 & 0.011 & 0.03 & $10^{-13}$ 
        & 1.41 & 1.42 & 0.042 & 0.06 & $10^{-13}$ 
        & 38.7 & 45.1 & 0.276 & 1.65 & $10^{-11}$ \\ \hline
      $50,000$    & 0.34 & 0.65 & 0.030 & 0.13 & $10^{-13}$ 
        & 3.31 & 3.43 & 0.082 & 0.15 & $10^{-12}$ 
        & 87.8 & 97.8 & 0.578 & 2.28 & $10^{-11}$ \\ \hline
      $100,000$   & 0.85 & 1.44 & 0.059 & 0.22 & $10^{-12}$ 
        & 6.54 & 6.95 & 0.177 & 0.31 & $10^{-12}$ 
        & 164  & 195  & 1.24  & 3.84 & $10^{-10}$ \\ \hline
      $200,000$   & 1.56 & 3.12 & 0.147 & 0.44 & $10^{-13}$ 
        & 14.1 & 15.9 & 0.395 & 0.59 & $10^{-11}$ &&&&&\\ \hline
      $500,000$   & 4.72 & 8.33 & 0.363 & 0.94 & $10^{-12}$ 
        & 38.2 & 42.1 & 1.12  & 1.69 & $10^{-11}$ &&&&&\\ \hline
      $1,000,000$ & 10.9 & 17.1 & 0.814 & 1.94 & $10^{-12}$ 
        & 79.9 & 90.3 & 2.38  & 3.39 & $10^{-11}$ &&&&&\\ \hline
      \end{tabular}
  }
    \end{center}
\end{table*}

\subsection{Exponential covariance}

Covariances functions of the form
\begin{equation}
k(r_i,r_j) = \exp(-|r_i-r_j|)
\end{equation}
are known as exponential covariance functions. One-dimensional numerical
results are presented in Table~\ref{tab_exp}.
 Figure~\ref{Figure_Rest} compares scaling for various kernels.

\begin{table}[!htbp]
    \caption{ \footnotesize Timings for one-dimensional exponential covariance
      functions. The matrix entry is given as $C_{ij} = \delta_{ij} +
      \exp(-|r_i-r_j|)$, where $r_i$ are random uniformly distributed
      points in the interval $[-3, 3]$.}
\label{tab_exp}
  \rowcolors{1}{gray!30}{white}
  \begin{center}
    \begin{tabular}{|c|cccc|c|}
      \hline & \multicolumn{4}{|c|}{Time taken in seconds} & \\ \hline
      $n$ & Assembly & Factor & Solve & Det. & Error \\
      \hline
      $10^4$         & 0.13 & 0.06 & 0.003 & 0.02 & $10^{-13}$ \\ \hline
      $2\times 10^4$ & 0.23 & 0.11 & 0.008 & 0.03 & $10^{-13}$ \\ \hline
      $5\times 10^4$ & 0.64 & 0.32 & 0.020 & 0.10 & $10^{-12}$ \\ \hline
      $10^5$         & 1.41 & 0.70 & 0.039 & 0.23 & $10^{-13}$ \\ \hline
      $2\times 10^5$ & 2.86 & 1.42 & 0.076 & 0.42 & $10^{-12}$ \\ \hline
      $5\times 10^5$ & 8.63 & 3.47 & 0.258 & 0.67 & $10^{-12}$ \\ \hline
      $10^6$         & 18.8 & 8.05 & 0.636 & 1.35 & $10^{-12}$ \\ \hline
    \end{tabular}
  \end{center}
\end{table}

\subsection{Inverse Multiquadric and Biharmonic}

\begin{figure}
  \begin{center}
    \includegraphics[width=.95\linewidth]
                    {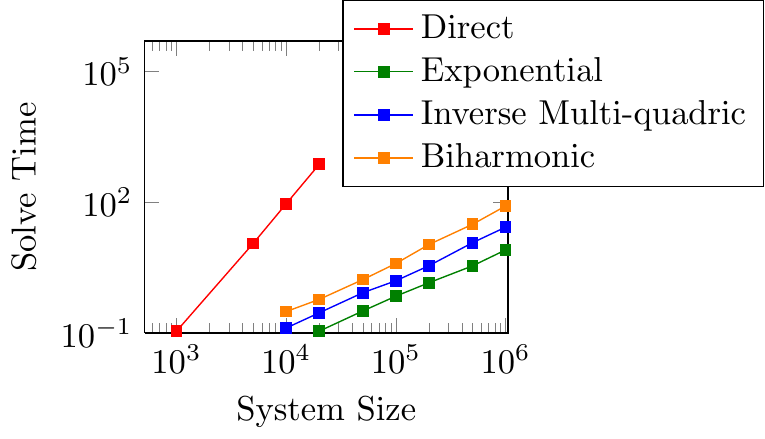}
  \end{center}
  \caption{ \footnotesize Comparison of time required to factorize a variety of
    covariance matrices arising from exponential, inverse
    multiquadric, and biharmonic covariance kernels for
    one-dimensional data.}
  \label{Figure_Rest}
\end{figure}

The inverse multiquadric and biharmonic kernel (also known as the thin
plane spline) are frequently used in radial basis function
interpolation and {\it kriging} in geostatistics. These kernels are
given by the formulae
\begin{equation}
\begin{aligned}
k_1(r_i,r_j) &= \frac{1}{\sqrt{1+| r_i - r_j |^2}}, \\
k_2(r_i,r_j) &= | r_i - r_j|^2 \log | r_i - r_j |,
\end{aligned}
\end{equation}
respectively. Timing results are presented in Tables~\ref{tab_inv}
and~\ref{tab_bi}, and comparison with the exponential kernel is
shown in Figure~\ref{Figure_Rest}.

\begin{table}[!htbp]
    \caption{ \footnotesize Timings for one-dimensional inverse multiquadric
      covariance functions. The matrix entry is given as $C_{ij} =
      \delta_{ij} + 1/\sqrt{1+{| r_i - r_j |^2}}$, where $r_i$ are
      random uniformly distributed points in the interval $[-3, 3]$.}
\label{tab_inv}
  \rowcolors{1}{gray!30}{white}
  \begin{center}
    \begin{tabular}{|c|cccc|c|}
      \hline & \multicolumn{4}{|c|}{Time taken in seconds} & \\ \hline
      $n$ & Assembly & Factor & Solve & Det. & Error \\
      \hline
      $10^4$         & 0.11 & 0.13 & 0.006 & 0.02 & $10^{-13}$ \\ \hline
      $2\times 10^4$ & 0.17 & 0.29 & 0.017 & 0.04 & $10^{-13}$ \\ \hline
      $5\times 10^4$ & 0.47 & 0.84 & 0.037 & 0.12 & $10^{-12}$ \\ \hline
      $10^5$         & 1.07 & 1.58 & 0.072 & 0.21 & $10^{-12}$ \\ \hline
      $2\times 10^5$ & 2.18 & 3.49 & 0.158 & 0.44 & $10^{-12}$ \\ \hline
      $5\times 10^5$ & 6.43 & 11.8 & 0.496 & 0.71 & $10^{-11}$ \\ \hline
      $10^6$         & 14.2 & 26.8 & 1.02  & 1.49 & $10^{-11}$ \\ \hline
    \end{tabular}
  \end{center}
\end{table}
%
%

\begin{table}[!htbp]
    \caption{ \footnotesize Timings for one-dimensional biharmonic covariance
      function. The matrix entry is given as $C_{ij} = 2\delta_{ij} +
      {| r_i - r_j|^2} \log (| r_i - r_j |)$, where $r_i$ are
      random uniformly distributed points in the interval $[-3, 3]$.}
\label{tab_bi}
  \rowcolors{1}{gray!30}{white}
  \begin{center}
    \begin{tabular}{|c|cccc|c|}
      \hline & \multicolumn{4}{|c|}{Time taken in seconds} & \\ \hline
      $n$ & Assembly & Factor & Solve & Det. & Error \\
      \hline
      $10^4$         & 0.28 & 0.31 & 0.015 & 0.03 & $10^{-12}$ \\ \hline
      $2\times 10^4$ & 0.61 & 0.59 & 0.028 & 0.06 & $10^{-12}$ \\ \hline
      $5\times 10^4$ & 1.62 & 1.68 & 0.067 & 0.15 & $10^{-12}$ \\ \hline
      $10^5$         & 3.61 & 3.93 & 0.123 & 0.34 & $10^{-12}$ \\ \hline
      $2\times 10^5$ & 8.03 & 10.7 & 0.236 & 0.65 & $10^{-12}$ \\ \hline
      $5\times 10^5$ & 26.7 & 31.2 & 0.632 & 1.41 & $10^{-11}$ \\ \hline
      $10^6$         & 51.3 & 81.9 & 1.28  & 3.40 & $10^{-11}$ \\ \hline
    \end{tabular}
  \end{center}
\end{table}

\subsection{Scaling in high dimensions}

In this section we report results on the scaling of the algorithm
described in this paper when the data (independent variables, $\bx$)
lie in high Euclidean dimensions. We perform two experiments. First,
we run our algorithm on data lying in the hypercube $[-3,3]^d$ for
various values of $d$. In this scenario, we actually see an {\it
  increase} in computational speed and accuracy as $d$ is increased
after some point. These results are reported in Table~\ref{tab_dim}.
However, this is not a fair result. As $d$ increases, the expected
value of $r = ||r_i - r_j||$ increases, causing, at least in the case
of an unscaled Gaussian covariance kernel, for many matrix entries to
be very close to zero.

The second experiment was with the same set of parameters, except the
data are located in the {\it scaled} hypercube,
$[-3/\sqrt{d},3/\sqrt{d}]^d$. These results are reported in
Table~\ref{tab_dim_scale}. This is equivalent to rescaling Euclidean
distance, or rescaling the covariance kernel.  We see that the
scalings for this experiment saturate once $d \approx 10$ due to the
fact that we are not increasing the number of data points $n$ along
with $d$. For fixed $n$, the data become very sparse as $d$ increases
and the ranks of the off-diagonal blocks in the associated covariance
matrix remain the same. What we mean to say by this is that $5000$
points in a ten-dimensional space is massively under-sampling any sort
of spatial structure, this is equivalent to grid of about two or three
points per dimension ($2.34^{10} \approx 5000$).  Running the
algorithm with this type of data is equivalent to doing dense linear
algebra since the ranks of all off-diagonal blocks are close to
full. This is a manifestation of the curse of dimensionality.
However, it's possible that one may encounter both types of data
(scaled vs. unscaled) in real-world situations.

These experiments were run on a faster laptop, namely a
MacBook Pro with a 3.0GHz Intel Core i7
processor and 16 GB 1600 MHz DDR3 RAM. No sophisticated software
optimizations were made.

\begin{table}[!htbp]
      \caption{ \footnotesize Timings for Gaussian covariance
        functions in $d$-dimensions. The matrix
        entries are given as $C_{ij} = 2\delta_{ij}
        + \exp \left(-|| r_i - r_j||^2\right)$, where $r_i$ are random
        uniformly distributed points in the unscaled interval $[-3,
        3]^d$
        (where $d = 1, \ldots, 64$).}
\label{tab_dim}
  \rowcolors{1}{gray!30}{white}
  \begin{center}
  \resizebox{\hsize}{!}{
    \begin{tabular}{|cc|cccc|c|}
      \hline & & \multicolumn{4}{|c|}{Time taken in seconds} & \\ \hline
      $n$ & Dim $d$ & Assembly & Factor & Solve & Det. & Error \\
      \hline
      $5000$ & 1  & 0.02 & 0.087 & 0.002 & 0.001 & $10^{-12}$ \\ \hline
      $5000$ & 2  & 1.17 & 2.043 & 0.027 & 0.014 & $10^{-11}$ \\ \hline
      $5000$ & 4  & 9.67 & 21.25 & 0.073 & 0.088 & $10^{-13}$ \\ \hline
      $5000$ & 8  & 10.4 & 7.773 & 0.041 & 0.039 & $10^{-08}$ \\ \hline
      $5000$ & 16 & 0.09 & 0.168 & 0.002 & 0.002 & $10^{-17}$ \\ \hline
      $5000$ & 32 & 0.15 & 0.061 & 0.004 & 0.001 & $10^{-16}$ \\ \hline
      $5000$ & 64 & 0.21 & 0.059 & 0.001 & 0.001 & $10^{-16}$ \\ \hline
    \end{tabular}
  }
\end{center}
\end{table}

\begin{table}[!htbp]
      \caption{ \footnotesize Timings for Gaussian covariance
        functions in $d$-dimensions. The matrix
        entries are given as $C_{ij} = 2\delta_{ij}
        + \exp \left(-|| r_i - r_j||^2\right)$, where $r_i$ are random
        uniformly distributed points in the scaled interval $[-3/\sqrt{d},
        3/\sqrt{d}]^d$
        (where $d = 1, \ldots, 64$).}
\label{tab_dim_scale}
  \rowcolors{1}{gray!30}{white}
  \begin{center}
  \resizebox{\hsize}{!}{
    \begin{tabular}{|cc|cccc|c|}
      \hline & & \multicolumn{4}{|c|}{Time taken in seconds} & \\ \hline
      $n$ & Dim $d$ & Assembly & Factor & Solve & Det. & Error \\
      \hline
      $5000$ & 1 &  0.02 & 0.087 & 0.002 & 0.001 & $10^{-12}$ \\ \hline
      $5000$ & 2 &  0.70 & 1.157 & 0.017 & 0.005 & $10^{-11}$ \\ \hline
      $5000$ & 4 &  33.6 & 108.2 & 0.165 & 0.360 & $10^{-12}$ \\ \hline
      $5000$ & 8 &  40.8 & 138.4 & 0.186 & 0.359 & $10^{-15}$ \\ \hline
      $5000$ & 16 & 41.9 & 138.7 & 0.182 & 0.358 & $10^{-15}$ \\ \hline
      $5000$ & 32 & 43.4 & 174.1 & 0.174 & 0.354 & $10^{-15}$ \\ \hline
      $5000$ & 64 & 42.2 & 142.2 & 0.181 & 0.346 & $10^{-15}$ \\ \hline
    \end{tabular}
  }
\end{center}
\end{table}

\subsection{Regression performance}

In this section we demonstrate the relationship between various
parameters in our algorithm and regression performance.  The two main
parameters that need to be set in our algorithm are the factorization
precision $\epsilon$ (see Figure~\ref{fig-algorithm}) and the maximum
size of the smallest sub-matrix on the finest level, $p_{max}$.  For a
fixed $\epsilon$, changing $p_{max}$ does {\em not} affect the $RMSE$
(root-mean-square error) of the regression (up to machine precision
errors), it merely affects the overall runtime. For sufficiently large
$p_{max}$, the scheme ceases to be a multi-level algorithm.  We merely
state this as a fact, and do not report the data.  We set $p_{max} =
20$ in the following numerical experiment.

\begin{figure}
  \begin{center}
    \includegraphics[width=.95\linewidth]{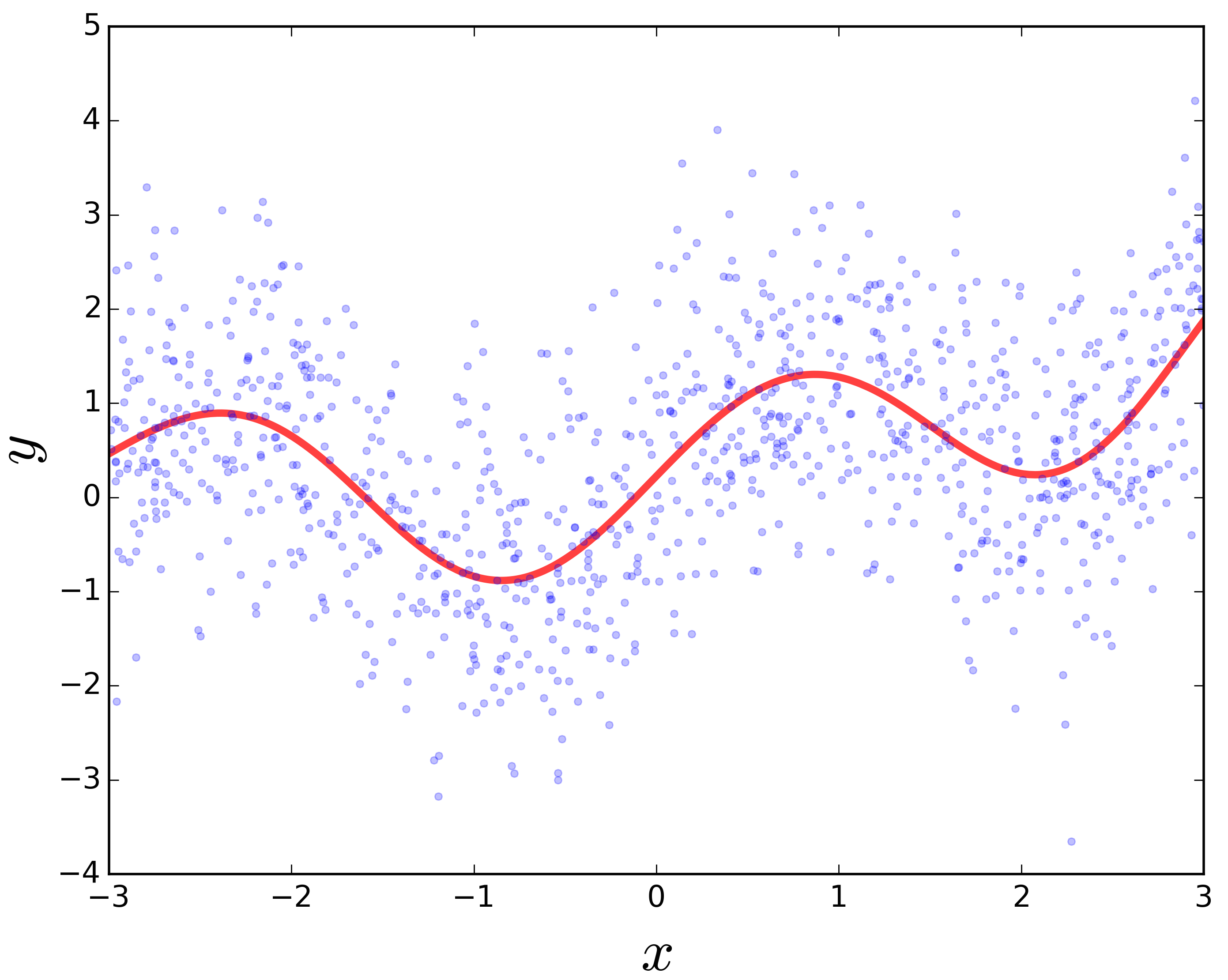}
  \end{center}
  \caption{ \footnotesize Regression obtained from Gaussian process
    prior for synthetic data drawn from the model in equation~\eqref{eq_model}.}
  \label{fig_regression}
\end{figure}

\begin{figure}
  \begin{center}
    \includegraphics[width=.95\linewidth]{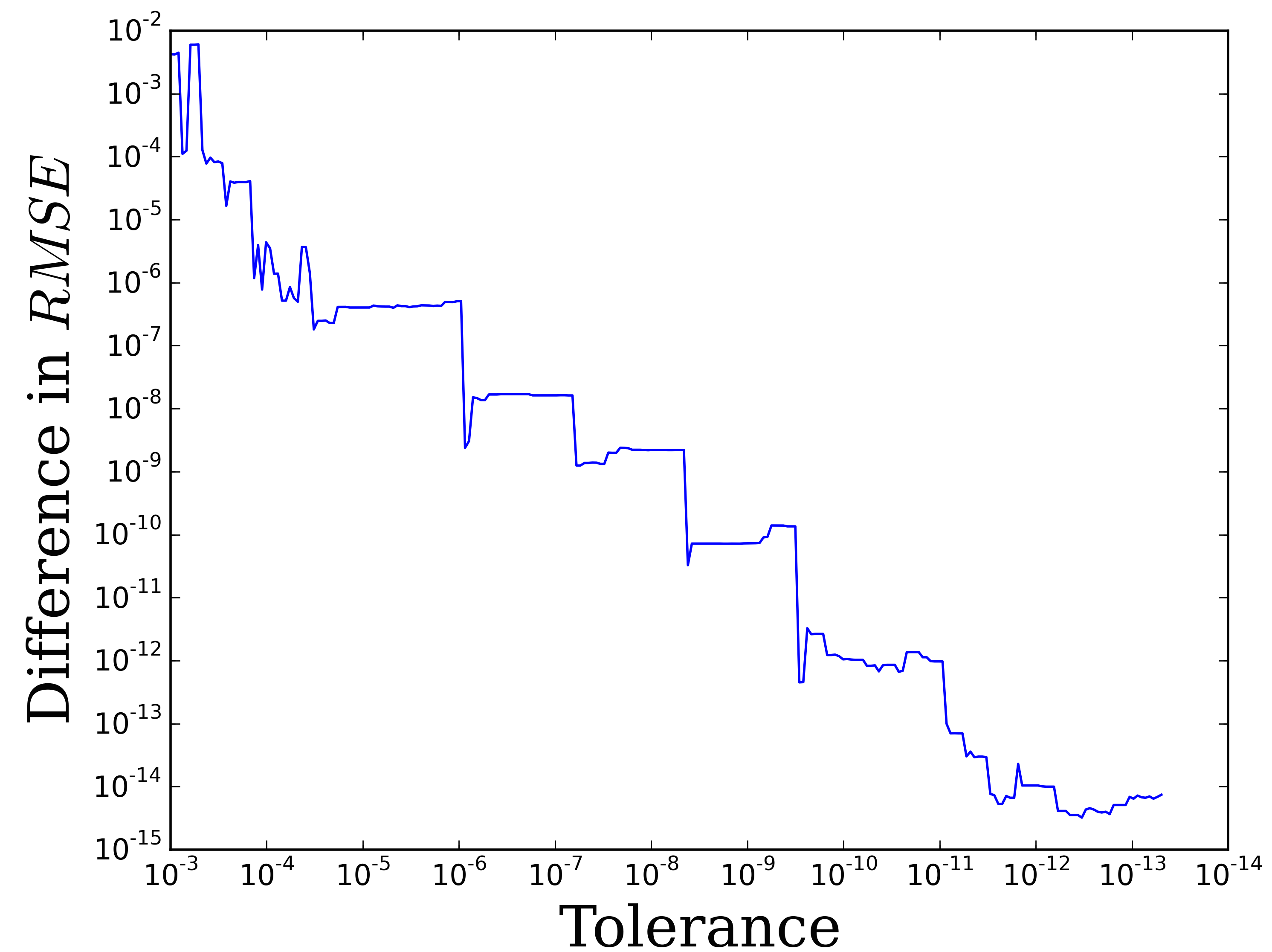}
  \end{center}
  \caption{ \footnotesize Difference in the root-mean-square error of
    the Gaussian process regression using the algorithm of this paper
    and that of the direct calculation as a function of the
    factorization precision $\epsilon$.}
  \label{fig_rmse}
\end{figure}

However, for varying values of $\epsilon$, we present the difference
between the $RMSE$ for the exact (dense linear algebra) regression and
the regression obtained using the matrix factorization algorithm of
this paper.  We generate $1024$ data-points $\{x_j, y_j\}$ from the
model:
\begin{equation}\label{eq_model}
y = \sin(2x) + \frac{1}{8} e^x + \epsilon,
\end{equation}
where $\epsilon \sim \mathcal N(0,\sigma^2_\epsilon)$ and $x_j$ is
chosen randomly in the interval $(-3,3)$. The non-parametric
regression curve (or estimate) under a Gaussian process prior (with
zero mean) is then calculated at the same points $x_j$ as
in~\eqref{eq_reg1} and~\eqref{eq_reg2}:
\begin{equation}
  \hat \by = K(\bx)
\left(\sigma^2_\epsilon I + K(\bx) \right)^{-1} \by,
\end{equation}
where we have chosen our covariance kernel to be consistent with
the previous numerical experiments:
\begin{equation}
  k(x,x') = e^{-(x-x')^2}.
\end{equation}
No effort was made to adapt the covariance kernel to the synthetic
data.  For $\sigma_\epsilon = 1.0$, the inferred curve through the
data is shown in Figure~\ref{fig_regression}. If dense linear algebra
is used to invert $I+K$ to obtain the {\em exact} estimate
$\by_{exact}$, the  $RMSE$ for this data is:
\begin{equation}
  \begin{aligned}
    RMSE &= \frac{\| \by - \hat \by_{exact} \|}{\sqrt{1024}} \\
    &= 0.9969571605921435.
  \end{aligned}
\end{equation}
Figure~\ref{fig_rmse} shows the {\em absolute difference} between the
$RMSE$ obtained from dense linear algebra and that obtained from our
accelerated scheme as a function of factorization precision
$\epsilon$. Unsurprisingly, from this plot we determine that, indeed,
the difference in regression performance (at least as measured by
$RMSE$) is proportional to $\epsilon$.

\section{Conclusions}
\label{sec-conclusions}

In this paper, we have presented a fast, accurate, and nearly optimal
hierarchical direct linear algebraic algorithms for computing
determinants, inverses, and matrix-vector products involving
covariances matrices encountered when using Gaussian
processes. Similar matrices appear in problems of classification and
prediction; our method carries over and applies equally well to these
problems.  Previous attempts at accelerating these calculations
(inversion and determinant calculation) relied on either sacrificing
fidelity in the covariance kernel (e.g. thresholding), constructing a
global low-rank approximation to the covariance kernel, or paying the
computational penalty of dealing with dense, full-rank covariance
matrices. Our HODLR-based algorithm obviates the need for this
compromise. Our {\it observation} that many covariance matrices of
mathematical statistics have fine-grained, compressible hierarchical
structure that provides access to the inverse may find use in 
many applications in the future.

The source code for the algorithm has been made available on
GitHub. The HODLR package for solving linear systems and computing
determinants is available at
\url{https://github.com/sivaramambikasaran/HODLR}~\cite{HODLR} and the
Python Gaussian process package~\cite{dfm_gp}, {\it george}, has been
made available at \url{https://github.com/dfm/george}.  Both packages
are open source, the HODLR package is released under the MPL2.0
license and {\it george} is released under the MIT license. Details on
using these packages are available at their respective online
repositories.

In its present form, our method degrades in performance when the
$n$-dimensional data has a covariance function based on points in
$\mathbb R^d$ with $d >3$, as well as when the covariance function is
oscillatory.
Part of the performance loss cannot be avoided due to the curse of
dimensionality. High-dimensional data is simply more complicated than
low-dimensional data causing the off-diagonal blocks to have larger
ranks (at least in the scenario of more and more data samples).
The other part of the performance loss is in the compression.
For high-dimensional data, analytic interpolatory low-rank
approximations will provide faster and more robust approximations.
Extensions of our approach to these cases is a subject of
current research. We are also investigating high-dimensional
anisotropic quadratures for marginalization and moment computation.

\ifCLASSOPTIONcompsoc
  \section*{Acknowledgments}
\else
  \section*{Acknowledgment}
\fi

The authors would like to thank Iain Murray for several useful
and detailed discussions.

\bibliographystyle{IEEEtranS}
\bibliography{./bibliography.bib}

\ifCLASSOPTIONcaptionsoff
  \newpage
\fi

\begin{IEEEbiographynophoto}{Sivaram Ambikasaran}
obtained his Bachelor's and Master's in Aerospace Engineering from
Indian Institute of Technology Madras in $2007$. Thereafter, he
received his Master's in Statistics, Master's \& Ph.D. in
Computational Mathematics from Stanford University in $2013$, where he
worked on {\it{fast direct solvers for large dense matrices}}. He is
now a Courant Instructor at New York University. His research focuses
on designing efficient, fast, scalable algorithms for mathematical
problems arising out of physical applications.
\end{IEEEbiographynophoto}

\begin{IEEEbiographynophoto}{Daniel Foreman-Mackey}
is a Ph.D.\ candidate in the Physics department at New York University. He
received a B.Sc.\ from McGill University in 2008 and a M.Sc.\ in physics from
Queen's University (Canada) in 2010. His research is focused on comprehensive
probabilistic data analysis projects in astrophysics, and specifically the
field of exoplanets.
\end{IEEEbiographynophoto}
\begin{IEEEbiographynophoto}{Leslie Greengard}
received a B.A. degree in Mathematics from Wesleyan University in
1979, a Ph.D. degree in Computer Science from Yale University in 1987,
and an M.D. degree from Yale University in 1987.  From 1987-1989 he
was an NSF Postdoctoral Fellow at Yale University and at the Courant
Institute of Mathematical Sciences, NYU, where he has been a faculty
member since 1989. He was the Director of the Courant Institute from
2006-2011 and is presently Director of the Simons Center for Data
Analysis at the Simons Foundation. His research interests include fast
algorithms, acoustics, electromagnetics, elasticity, heat transfer,
fluid dynamics, computational biology and medical imaging.
Prof. Greengard is a member of the National Academy of Sciences and
the National Academy of Engineering.
\end{IEEEbiographynophoto}
\begin{IEEEbiographynophoto}{David W. Hogg}
obtained an S.B. from the Massachusetts Institute of Technology and a
PhD in physics from the California Institute of Technology.  His
research has touched on large-scale structure in the Universe, the
structure and formation of the Milky Way and other galaxies, and
extra-solar planets.  He currently works on engineering and data
analysis aspects of large astronomical projects, with the goal of
increasing discovery space and improving the precision and efficiency
of astronomical measurements.
\end{IEEEbiographynophoto}
\begin{IEEEbiographynophoto}{Michael O'Neil}
received an A.B. in Mathematics from Cornell University in 2003, and a
Ph.D. in Applied Mathematics from Yale University in 2007 where he
developed the generalized framework for {\it butterfly-type}
algorithms, widely applicable to oscillatory integral operators. He
was a postdoctoral fellow at the Courant Institute of Mathematical
Sciences, NYU from 2010-2014. Presently, he is Assistant Professor of
Mathematics at the Courant Institute, NYU and the Polytechnic School
of Engineering, NYU. His
research focuses on problems in computational electromagnetics,
acoustics, and magnetohydrodynamics, integral equations, fast
analysis-based algorithms, and computational statistics.
\end{IEEEbiographynophoto}





\end{document}